\documentclass[11pt]{article}
\usepackage{latexsym}
\usepackage{epsfig,enumerate, url}
\usepackage{amsmath,amsthm,amssymb, bbm}
\usepackage{algorithm2e}
\usepackage{appendix}
\usepackage{hyperref}
\usepackage{comment}
\usepackage{float}
\usepackage{enumitem}
\usepackage{caption}

\usepackage[usenames]{color}\usepackage{graphicx}
\numberwithin{equation}{section}
\definecolor{brown}{cmyk}{0, 0.72, 1, 0.45}
\definecolor{grey}{gray}{0.5}

\renewcommand{\epsilon}{\varepsilon}

\newtheorem{opq}{Open Question}
\newtheorem*{conjecture*}{Conjecture}
\newtheorem{theorem}{Theorem}[section]
\newtheorem*{theorem*}{Theorem}

\newcommand{\rbrac}[1]{\left(#1\right)}
\newcommand{\cbrac}[1]{\left\{#1\right\}}
\newcommand{\sbrac}[1]{\left[{ #1}\right]}

\newcommand{\bfo}{\mathbbm{1}}

\def\sm{\setminus}

\def\p{\propto}
\def\t{\tau}
\def\l{\ell}

\def\op{\mbox{op}}
\def\Op{\mbox{Op}}

\def\uar{u.a.r.}
\newcommand{\nn}{\nonumber}

\def\E{\mathbb{E}}

\def\yb{{\bf {y}}}

\def\zb{{\bf {z}}}

 \newcommand{\tbf}[1]{\textbf{#1}}

\newcommand{\mc}[1]{\mathcal{#1}}

\title{Larger matchings and independent sets in regular uniform hypergraphs of high girth}
\author{Deepak Bal\thanks{Department of Mathematics, Montclair State University, Montclair, NJ. {\tt deepak.bal@montclair.edu}} \and Patrick Bennett\thanks{Department of Mathematics, Western Michigan University, Kalamazoo, MI. {\tt patrick.bennett@wmich.edu} \\ Research supported in part by Simons Foundation Grant \#426894.}}
\date{}

\begin{document}

\maketitle

\begin{abstract}
    In this note we analyze two algorithms, one for producing a matching and one for an independent set, on $k$-uniform $d$-regular hypergraphs of large girth. As a result we obtain new lower bounds on the size of a maximum matching or independent set in such hypergraphs. 
\end{abstract}

\section{Introduction}

A \tbf{hypergraph} is a pair $\mc{H} =(V, E)$ with $E \subseteq 2^V$.  We call the elements of $V$ {\bf vertices} and elements of $E$ \tbf{edges}. We write $V(\mc{H})=V$ and $E(\mc{H})=E$. We say  $\mc{H}$ is $k$-uniform if each edge has exactly $k$ vertices. We say $\mc{H}$ is $d$-regular if each vertex is in $d$ edges. We let $\mc{H}_k(n, d)$ be the random $d$-regular $k$-uniform hypergraph on vertex set $V=[n]$. In other words, $\mc{H}_k(n, d)$ is chosen uniformly at random (u.a.r.) from all $d$-regular $k$-uniform hypergraphs on $V$. In this note we are interested in fixing $k, d$ and letting $n$ go to infinity. We will always assume $k$ divides $dn$, which for large $n$ is sufficient to guarantee the existence of $d$-regular $k$-uniform hypergraphs on $n$ vertices.

A \tbf{matching} in a hypergraph $\mc{H}$ is a set of pairwise-disjoint edges. An \tbf{independent set} in $\mc{H}$ is a set of vertices that does not contain any edges.  There are many ways to define cycles in hypergraphs, but for us, a \tbf{cycle of length $\l$} in $\mc{H}$ is a sequence $v_1, e_1, v_2, e_2, \ldots, v_{\l}, e_{\l}, v_{\l+1}=v_1$ that alternates vertices and edges such that  all vertices and edges in the sequence are distinct (with the exception that $v_{\l+1}=v_1$), and  for each $j$ the edge $e_j$ contains vertices $v_j$ and $v_{j+1}$. The cycles we define here are often called Berge cycles. The \tbf{girth} of a hypergraph $\mc{H}$ is the minimum length of a cycle (or $\infty$ if there are no cycles). Note that for $d, k \ge 2$ any $k$-uniform $d$-regular hypergraph $\mc{H}$ has a cycle and thus finite girth.

In this note we analyze two algorithms, one for producing a matching and one for an independent set, on uniform, regularhypergraphs of large girth. They are both greedy algorithms, in the sense that each algorithm builds its output one element at a time (i.e. we build a matching one edge a time, or an independent set one vertex at a time) and we never consider the possibility of removing an element once we have included it. At each step we must avoid making a choice that conflicts with the previous choices to which we committed. We therefore make some attempt at each step to maximize the number of choices we will have in the future. In a greedy matching algorithm, each time we choose a new edge $e$ to go in our matching we are forced to eliminate every edge that meets $e$ from future consideration. Thus it makes sense to try to choose $e$ so that it does not meet too many other edges. Similarly we take vertex degrees into account for our independent set algorithm. Our algorithms are stated formally after this paragraph. Both algorithms above use the standard operation of {\bf vertex deletion}. When a vertex is deleted from a hypergraph, all edges containing that vertex are also deleted. 

\begin{center}
\begin{minipage}{.45\textwidth}
\begin{algorithm}[H]\label{algo:match}
  \SetKwInOut{Input}{Input}
   \SetKwInOut{Output}{Output}
	\caption{\textsc{degree-greedy matching process}}
	\Input{Hypergraph $\mc{H}=(V,E)$}
	\Output{Matching $M$}
	$M = \emptyset$\;
	\While{$E\neq \emptyset$}
	{Select $v \in V$ \uar~from all vertices of minimum positive degree\;
 Select an edge $e \ni v$ \uar\;
	 $M	\gets M \cup \cbrac{e}$\;
    $V \gets V \sm e$\;
	}
	\Return $M$\;
\end{algorithm} 
\end{minipage}
    \begin{minipage}{.45\textwidth}
\begin{algorithm}[H]\label{algo:ind}
  \SetKwInOut{Input}{Input}
   \SetKwInOut{Output}{Output}
	\caption{\textsc{degree-greedy independent process}}
	\Input{Hypergraph $\mc{H}=(V,E)$}
	\Output{Independent set $I$}
	$I = \emptyset$\;
	\While{$V \neq \emptyset$}
	{Select $v \in V$ \uar~from all vertices of minimum degree\;
	 $I	\gets I \cup \cbrac{v}$\;
	 $V \gets V\setminus \cbrac{v'\in V : \exists e\in E, e\subseteq I \cup \{v'\}}$\;
	}
	\Return $I$\;
\end{algorithm} 
\end{minipage}
\end{center}

These algorithms are {\em degree-greedy} versions of processes that the present authors analyzed in \cite{BB}, which we call the {\bf matching process} and the {\bf independent process}. In the matching process (resp. independent process), at each step we choose an edge (resp. vertex) uniformly at random from all edges (resp. vertices) whose addition would still be a matching (resp. independent set). 

Both Algorithms \ref{algo:match} and \ref{algo:ind}  are inspired by algorithms which have been analyzed on graphs. Karp and Sipser \cite{KS} analyzed a matching algorithm which adds a pendant edge whenever one exists, or else adds a uniformly random edge. Their analysis was performed on $G(n,p)$ where $p=c/n$. This algorithm was further analyzed by Aaronson, Frieze and Pittel \cite{AFP} on $G(n,p)$ and subsequently analyzed on random graphs with a fixed degree sequence by Bohman and Frieze \cite{BF11}. In particlar, for random regular graphs, this algorithm produces an {\bf almost-perfect matching}, i.e. a matching that covers all but $o(n)$ vertices.  A fully degree prioritized algorithm, MINGREEDY, was analyzed by Frieze, Radcliffe and Suen \cite{FRS} on random 3-regular graphs and shown to produce a matching covering all but at most $\tilde{\Theta}(n^{1/5})$ vertices. Algorithm \ref{algo:match} is the natural extension of the MINGREEDY algorithm to hypergraphs. 

For independent sets, a degree-greedy algorithm was analyzed on random 3-regular graphs by Frieze and Suen \cite{FS}. This algorithm produced an independent set of size  
at least $(6\log(3/2) - 2 - \epsilon) n$. Independently, Wormald \cite{Wor95} analyzed a degree-greedy algorithm for random $r$-regular graphs and provided numerical lower bounds on the independent sets produced by this algorithm.  Algorithm \ref{algo:ind} is the natural extension to hypergraphs.

In this area, the so-called differential equation method of Wormald \cite{Wor95} has been extremely useful. The method allows one to prove that a system of random variables is {\bf dynamically concentrated} around their expected trajectories throughout the course of a random process. The term ``dynamic concentration'' refers to the fact that these random variables change over time, and at each step they are concentrated around their trajectories (which obviously also change over time). The trajectories are given by the solution to a system of differential equations, hence the name.

Perhaps more than random regular graphs, there has been a lot of interest in properties of arbitrary regular graphs of high girth. Graphs of girth $g$ are ``locally treelike'' in the sense that exploring the first $ g/2 -1$ neighborhoods at any vertex yields a tree. Furthermore if the graph is regular, these trees are all isomorphic. Thus, in a sense, all high-girth $d$-regular graphs look the same locally. Pippenger \cite{Pip} and independently Lauer and Wormald \cite{LW} analyzed the independent process on high-girth regular graphs, finding the expected size of the final independent set (Pippenger \cite{Pip} also did something similar for the matching process). Gamarnik and Goldberg \cite{GG10} established concentration of the size of the final independent set or matching produced by the processes (again working on regular graphs of high girth). Nie and Verstra\"ete \cite{NV} then analyzed the independent process on regular hypergraphs of large girth. The present authors \cite{BB} also gave lower bounds on the size of a maximum matching in high-girth regular hypergraphs, using a beautiful result of Krivelevich, M\'esz\'aros, Michaeli and Shikhelman \cite{KMMS}. The result in \cite{KMMS} essentially says that whenever two graphs ``look locally almost the same'' (more formally, the probability distribution given by exploring several neighborhoods of a random vertex is asymptotically the same for both graphs), the independent process performs almost the same in each graph.

Since random regular graphs have few short cycles, they ``look locally almost the same'' as high-girth regular graphs, and therefore one might optimistically hope that the performance of more general algorithms on high-girth regular graphs to be about the same as random regular graphs. The amazing result of Hoppen and Wormald \cite{HWJCTB, HWComb}, allows for such a transfer. These results of Hoppen and Wormald immediately gave many applications of the following form: first one analyzes the performance of an algorithm on random regular graphs (to give some bound on a graph parameter), and then applies their result to get the same bound for all high-girth regular graphs. Since then, several similar applications have been given. For example, this method was used by Cs\'oka \cite{C} to give the best known lower bound for the independence number of high-girth cubic graphs. 

To summarize, lower bounds have been attained for matchings and independent sets in high girth \emph{graphs} using greedy algorithms (Gamarnik and Goldberg \cite{GG10})  and then later degree prioritized algorithms (Hoppen and Wormald \cite{HWJCTB, HWComb}). Lower bounds have been attained (by Nie and Verstra\"{e}te \cite{NV}) for high girth \emph{hypergraphs} by using the random greedy algorithm. In this note, we aim to fill in the missing piece of this story: we analyze degree prioritized algorithms on random regular hypergraphs,  then we observe how the results of Hoppen and Wormald \cite{HWJCTB, HWComb} apply in the hypergraph setting. As a result we improve the best known lower bounds for independent sets and matchings in high-girth regular hypergraphs.

\begin{theorem}
    Let $\nu_{k,d}$ and $\alpha_{k,d}$ be as in Tables \ref{tbl1} and \ref{tbl4}.  W.h.p., $\alpha(\mathcal{H}_k(n,d)) \ge \alpha_{k,d}n$ and  $\nu(\mathcal{H}_k(n,d)) \ge \nu_{k,d}n$. Furthermore, if $\mc{H}$ is any $k$-uniform $d$-regular graph with $n$ vertices and sufficiently large girth, then we likewise have $\alpha(\mathcal{H}) \ge \alpha_{k,d}n$ and  $\nu(\mathcal{H}) \ge \nu_{k,d}n$.
\end{theorem}

The entries in Table \ref{tbl1} are an improvement on the previous best known lower bounds, which were proved by the present authors in \cite{BB} and are summarized in Table \ref{tbl2}. Table \ref{tbl3} gives the corresponding upper bounds. More precisely, each entry $b_{d, k}$ in Table \ref{tbl3} indicates that w.h.p.~any matching in $\mc{H}_k(n, d)$ must have at most $b_{d, k}n$ edges and therefore leave at least $(1-kb_{d, k})n$ vertices unmatched. Note that this implies there exist $k$-uniform $d$-regular hypergraphs of arbitrarily high girth such that every matching leaves $(1-kb_{d, k})n$ vertices unmatched. Likewise, Table \ref{tbl4} improves the lower bounds given by and present authors \cite{BB} (for random regular hypergraphs) and Nie and Verstra\"ete \cite{NV} (for high-girth regular hypergraphs). These previous lower bounds for independent sets are in Table \ref{tbl5}. The upper bounds in Table \ref{tbl6} were proved by the present authors in \cite{BB}. Note that in many of the cases listed in Tables \ref{tbl1}--\ref{tbl6}, the new lower bound is closer to the upper bound than it is to the previous lower bound.
\begin{table}[hbt!]
\centering
\begin{tabular}{ |c|c|c|c|c| } 
 \hline
  & $d=2$ & $d=3$ & $d=4$ & $d=5$ \\ 
 \hline
 $k=3$ & 0.274 & 0.284 & 0.291 & 0.296 \\ 
 \hline
 $k=4$ & 0.179 & 0.181 & 0.186 &0.190\\ 
 \hline
 $k=5$ & 0.128 & 0.127 & 0.130 & 0.132\\ 
 \hline
\end{tabular}
\caption{Values for $\nu_{k, d}$}\label{tbl1}
\end{table}

\begin{minipage}{.45\textwidth}
\centering
\begin{tabular}{ |c|c|c|c|c| } 
 \hline
  & $d=2$ & $d=3$ & $d=4$ & $d=5$ \\ 
 \hline
 $k=3$ & 0.250 & 0.250 &0.253 & 0.257\\ 
 \hline
 $k=4$ & 0.166 & 0.164 & 0.166 & 0.169 \\ 
 \hline
 $k=5$ & 0.120 & 0.117 & 0.118& 0.120\\ 
 \hline
\end{tabular}
\captionof{table}{Previous lower bounds for $\nu(\mc{H})$}\label{tbl2}
\end{minipage}
\begin{minipage}{.45\textwidth}
\centering
\begin{tabular}{ |c|c|c|c|c| } 
 \hline
  & $d=2$ & $d=3$ & $d=4$ & $d=5$ \\ 
 \hline
 $k=3$ & 0.307 & 0.316 & 0.324 & 0.330 \\ 
 \hline
 $k=4$ & 0.211 & 0.216 & 0.221 & 0.226 \\ 
 \hline
 $k=5$ & 0.156 & 0.158 & 0.162 & 0.166 \\ 
 \hline
\end{tabular}
 \captionof{table}{Upper bounds for $\nu(\mc{H})$}\label{tbl3}
\end{minipage}

\begin{table}[hbt!]
\centering
\begin{tabular}{ |c|c|c|c|c| } 
 \hline
  & $d=2$ & $d=3$ & $d=4$ & $d=5$ \\ 
 \hline
 $k=3$ & 0.666 & 0.626 & 0.600 & 0.564\\ 
 \hline
 $k=4$ & 0.749 & 0.720 & 0.694 & 0.672\\ 
 \hline
 $k=5$ & 0.799 & 0.777 & 0.755& 0.737\\ 
 \hline
\end{tabular}
\caption{Values for $\alpha_{k, d}$}\label{tbl4}
\end{table}

\begin{minipage}{.45\textwidth}
\centering
\begin{tabular}{ |c|c|c|c|c| } 
 \hline
  & $d=2$ & $d=3$ & $d=4$ & $d=5$ \\ 
 \hline
 $k=3$ & 0.614 & 0.567 & 0.531 & 0.503 \\ 
 \hline
 $k=4$ & 0.708 & 0.670 & 0.640 & 0.616\\ 
 \hline
 $k=5$ & 0.765 & 0.733 & 0.708 & 0.688 \\ 
 \hline
\end{tabular}
 \captionof{table}{Previous lower bounds for $\alpha(\mc{H})$}\label{tbl5}
\end{minipage}
\begin{minipage}{.45\textwidth}
\centering
\begin{tabular}{ |c|c|c|c|c| } 
 \hline
  & $d=2$ & $d=3$ & $d=4$ & $d=5$ \\ 
 \hline
 $k=3$ & 0.667 & 0.651 & 0.624 & 0.600 \\ 
 \hline
 $k=4$ & 0.750 & 0.744 & 0.724 & 0.706\\ 
 \hline
 $k=5$ & 0.800 & 0.798 & 0.784 & 0.769\\ 
 \hline
\end{tabular}
 \captionof{table}{Upper bounds for $\alpha(\mc{H})$}\label{tbl6}
\end{minipage}

Although Tables \ref{tbl1}--\ref{tbl6} only apply when $ k \le 5$, $ d \le 5$, we provide a general method which can produce the appropriate bounds for Tables \ref{tbl1} and \ref{tbl4} for all fixed $k \ge 3, d \ge 2$. We have also uploaded the Maple files which we used to produce these values in the ancillary files on the arXiv version of this paper. The caveat is that for large $k, d$ this starts to become computationally expensive since it requires numerically solving larger and larger systems of differential equations. 

\section{Setting up the analysis for the processes} \label{sec:setup}

 In this section we describe the analysis of each process on the random regular hypergraph $\mc{H}_k(n, d)$ using the differential equation method. 

We generate our random hypergraph $\mc{H}_k(n, d)$ as follows. We will have two disjoint sets $P_A$ and $P_B$, each consisting of $nd$ \tbf{points}. $P_A$  is partitioned into $n$ parts of $d$ points each (the parts of this partition are called \tbf{vertices}). $P_B$ is partitioned into $nd/k$ parts of $k$ points each (these parts will be called  \tbf{edges}).   We generate a uniform random \tbf{pairing} $\mc{P}$ (a partition of $P_A \cup P_B$ into sets of size 2 each containing one point from $P_A$ and one from $P_B$).  We interpret the pairing as a hypergraph as follows: a vertex $v$ is contained in an edge $e$ if and only if there is a point in vertex $v$ that is paired with a point in the edge $e$.

We would like to restrict our attention to hypergraphs without ``loops" (i.e. edges that contain a vertex multiple times) or ``multi-edges" (i.e. two edges consisting of the same set of vertices). We call such hypergraphs \tbf{simple}. The main theorem in Blinovsky and Greenhill's paper \cite{BG16} can help us here. For fixed $k, d$ their main result gives
\[
P[\mc{H}_k(n, d) \mbox{ is simple}] = \Omega(1).
\]

At step $i$, say there are $Y_j = Y_j(i)$ many vertices with exactly $j$ unpaired points for $j=1, \ldots, d $. Say there are $Z_j=Z_j(i)$ edges with exactly $j$ unpaired points for $j=1, \ldots, k$. In a slight abuse of notation we will let $E=E(i)$ be the number of edges in the set $E$ at step $i$ (we will often use the same notation for a set and its cardinality). 
Let 
\begin{align}
 L_1& =\sum_{j=1}^d jY_j, \quad L_2 = \sum_{j=1}^d j(j-1)Y_j, \label{eqn:L}\\
 M_1& =\sum_{j=1}^{k} jZ_j, \quad M_2 = \sum_{j=1}^{k} j(j-1)Z_j.\label{eqn:M}
\end{align}
Since the number of unpaired points in $P_A$ is always the same as in $P_B$ we have $L_1 = M_1$.

\section{Analysis for the degree-greedy matching process}


We analyze Algorithm \ref{algo:match} on $\mc{H}_k(n, d)$. We start at step $0$ with the empty matching $M(0):=\emptyset$ and the hypergraph $\mc{H}(0):=\mc{H}_k(n, d)$. At each step $i \ge 1$ we do the following:
\begin{enumerate}
    \item choose a vertex $v_i$ uniformly at random from all vertices of smallest possible positive degree in $\mc{H}(i-1)$ (or if all vertices have degree 0 we terminate), 
    \item choose an edge $e_i$ uniformly at random from all edges containing $v_i$ and put $M(i):=M(i-1) \cup \{e_i\}$,
    \item form $\mc{H}(i)$ from $\mc{H}(i-1)$ by deleting all vertices in $e$.
\end{enumerate}

We now describe the above algorithm in terms of the pairing $\mc{P}$. We will only reveal parts of the pairing as the algorithm progresses. We start at step 0 with empty matching $M(0):=\emptyset$ and partial pairing $\mathcal{P}(0):=\emptyset$. By a partial pairing we mean a collection of disjoint pairs each containing one point from each of $P_A$ and $P_B$. At each step $i\ge 1$ we do the following:

\begin{enumerate}
    \item[(\tbf{M1})] From among all vertices with the smallest positive number $\delta$ of unpaired points, we choose one uniformly at random, say $v_i$. Pair all unpaired points of $v_i$.
    \item[(\tbf{M2})] Choose an edge $e_i$ uniformly at random from all edges containing a point that was just paired with a point from $v_i$. Set $M(i):=M(i-1) \cup \{e_i\}$.
    \item[(\tbf{M3})] Pair all remaining unpaired points in $e_i$. 
    \item[(\tbf{M4})] For any vertex $v'$ containing a point that was just paired 
    in (\tbf{M3}),
    we pair the remaining points in $v'$.
    \item[(\tbf{M5})] For any edge $e'$ containing a point that was just paired 
    in (\tbf{M4}) or in (\tbf{M1}), we pair the remaining points in $e'$. 
\end{enumerate}

Note that at the beginning of each step $i$ of this algorithm, each edge has either $0$ or $r$ unpaired points. Hence there is no need to track the $Z_j$'s in this algorithm. The expected one step changes in the $Y_j$'s are give by the following. Given a partial pairing $\mc{P}_i$ with $\delta = r$ we have

\begin{align}
      &\E\sbrac{Y_j(i+1)-Y_j(i) \mid  \mc{P}_i}  \nonumber\\
& = -\bfo_{j=r} - \frac{(k-1)jY_j}{L_1} +  \rbrac{r-1 +(k-1)\frac{L_2}{L_1}} (k-1) \rbrac{\frac{(j+1)Y_{j+1} - jY_j}{L_1}}+O\rbrac{\frac{1}{L_1}} \label{eq:mExpOneStepMatch}
\end{align}
Note that in one step (i.e. completing substeps (\tbf{M1})-(\tbf{M5})) we pair $O(1)$ points. Thus the random variables $Y_j, L_1, L_2$ change by at most $O(1)$. Thus, for example, when we pair a point in $P_B$ during step $i$, the probability that it is paired with a point in $Y_j$ is $\frac{jY_j(i) +O(1)}{L_1(i) + O(1)} = \frac{jY_j(i) }{L_1(i) }+O\rbrac{\frac{1}{L_1(i)}}$. For the remainder of this paragraph, we ignore such error terms since they are already accounted for by the term $O\rbrac{\frac{1}{L_1}}$.
The first term in \eqref{eq:mExpOneStepMatch} accounts for the loss of a vertex in $Y_j$ if the vertex chosen in (\tbf{M1}) is of degree $j$ (which happens if $j=r$).
The second term accounts for the expected number of vertices of degree $j$ paired to in (\tbf{M3}), keeping in mind that edge $e_i$ has $k-1$ unpaired points in (\tbf{M3}). 
The first factor of the third term describes the expected number of edges $e'$ referred to in (\tbf{M5}): $(r - 1)$ such $e'$ are paired to in (\tbf{M1}) and $(k-1)\sum_{\l=1}^{d}\frac{\l(\l-1)Y_\l}{L_1}=(k-1)\frac{L_2}{L_1}$ are paired to in (\tbf{M4}). Each such $e'$ has $k-1$ remaining unpaired points and the third factor of the third term is the expected change in $Y_j$ caused by these pairings.

Motivated by \eqref{eq:mExpOneStepMatch}, we define the following function. 
\begin{equation}
    f_{j, r}(\yb):= -\bfo_{j=r} - \frac{(k-1)jy_j}{\l_1} +  \rbrac{r-1 +(k-1)\frac{\l_2}{\l_1}} (k-1) \rbrac{\frac{(j+1)y_{j+1} - jy_j}{\l_1}},
\end{equation}
where
\begin{align}
 \l_1& =\sum_{i=1}^d iy_i, \qquad \l_2 = \sum_{i=1}^d i(i-1)y_i. \label{eqn:l}
\end{align}

 Of course the minimum degree $\delta$ of $\mc{H}$ is $\delta = d$ for the first step, and for the next several steps we expect to have $\delta = d-1$. After not too long we expect to see our first vertex of degree $d-2$, but for a while the algorithm should process these vertices at a rate that keeps them at or near zero. We (roughly) refer to these steps as \tbf{phase 1} of the algorithm, and we say that \tbf{phase 2} begins when the vertices of degree $d-2$ start to accumulate faster than the algorithm can process them. Indeed, typically this algorithm will have several \tbf{phases}, where phase $p$ consists of a large number of consecutive steps during which the minimum degree is almost always $d-p$ or $d-p-1$ and the number of vertices of degree $d-p-1$ stays relatively small. In order to keep the number of vertices of degree $d-p-1$ small during phase $p$, the algorithm must spend a significant proportion of steps processing them (and the rest of the steps processing vertices of degree $d-p$). Awkwardly, the expected one-step change in our variables \ref{eq:mExpOneStepMatch} depends on the degree of the vertex we process in that step. The standard way to handle this (see, e.g. Wormald \cite{WPhases}) is to determine the proportion of each type of step, and then write a ``blended'' version of our one-step change which is a sum of the possible one-step changes weighted according to the proportion of each type of step. We will let
 \begin{equation}
        \alpha_p(\yb) := f_{d-p-1, d-p}(\yb), \qquad \tau_p(\yb) := -f_{d-p-1, d-p-1}(\yb).
\end{equation}
In order to hold the vertices of degree $d-p-1$ near zero we must process them at the rate they appear. Thus the proportion of steps where we process a vertex of degree $d-p-1$ in phase $p$ should be about $\frac{\alpha_p}{\tau_p+\alpha_p}$, and of course then about $\frac{\tau_p}{\tau_p+\alpha_p}$ proportion of steps are spent processing vertices of degree $d-p$. The one exception is phase $d-1$ where we only process vertices of degree $1$ and never degree $0$ (recall we always process the vertex of minimum {\em positive} degree). 
The blended differential equations for phase $p$ are then
\begin{align}
    \frac{dy_j}{dx}& = F( \yb, j, p) := \begin{cases}
    \frac{\tau_p}{\tau_p+\alpha_p} f_{j, d-p}( \yb) + \frac{\alpha_p}{\tau_p+\alpha_p} f_{j, d-p-1}( \yb),& \;\;\; p \le d-2\\
    f_{j, 1}( \yb),& \;\;\; p =d-1.
    \end{cases}\label{eqn:mixedDEQY}
\end{align}

By now, this method of blending one-step changes for different types of steps during different phases of degree-greedy algorithms is standard. For random regular graphs this was formalized in a very general setting by Wormald \cite{WPhases}. Phase $1$ starts at $t_1=0$. For $p \ge 1$, phase $p$ ends when the algorithm is barely able to process the vertices of degree $d-p-1$ at the rate they accumulate, i.e. when $\tau_p(\yb) = -f_{d-p-1, d-p-1}(\yb)=0$. At the time that occurs, say $t_p$, we begin phase $p+1$ which then runs until time $t_{p+1}$ and so on. When phase $p$ begins we use the differential equation \eqref{eqn:mixedDEQY} together with initial conditions given by the state of the process at time $t_{p}$. The process ends when there are no more edges, i.e. when $\yb$ is the zero vector. We produced the values in Table \ref{tbl1} using Maple's numerical differential equation solver.

\section{Analysis for the degree-greedy independent process}

We analyze Algorithm \ref{algo:ind} on $\mc{H}_k(n, d)$. We start at step 0 with an empty independent set $\mc{I}(0):=\emptyset$ and the hypergraph $\mc{H}(0):=\mc{H}_k(n, d)$. At each step $i\ge 1$, we do the following:
\begin{enumerate}
    \item choose a vertex $v_i$ uniformly at random from all vertices of smallest degree in $\mc{H}(i-1)$. Put $\mc{I}(i) := \mc{I}(i-1)\cup \{v_i\}$.
    \item Remove $v_i$ from any edge containing $v_i$. In other words, any edge $e \ni v_i$ is replaced with the edge $e \sm \{v_i\}$.
    \item Some edges may now have only one vertex. For each such edge $e=\{w\}$ we say that $w$ is now \tbf{closed}. We delete all closed vertices and $v_i$ to form $\mc{H}_{i}$ (when we delete a vertex we also delete all incident edges).
\end{enumerate}

We now describe the above algorithm in terms of the pairing $\mc{P}$. We will only reveal parts of the pairing as the algorithm progresses. We start at step 0 with empty independent set $I(0):=\emptyset$. At each step $i\ge 1$ we do the following:
\begin{enumerate}
    \item[(\tbf{I1})] Choose a vertex $v_i$ uniformly at random from all vertices with the smallest number of unpaired points. Set $I(i):=I(i-1) \cup \{v_i\}$. Pair all unpaired points of $v_i$. 
    \item[(\tbf{I2})] For any edge $e$ containing a point that was just paired to a point from $v_i$, if $e$ now has only one unpaired point then we pair that point. 
    \item[(\tbf{I3})] For any vertex $w$ containing a point that was just paired in (\tbf{I2}), we pair the remaining points of $w$. We say that $w$ is closed now. 
    \item[(\tbf{I4})] For any edge $e'$ containing a point that was just paired in (\tbf{I3}), we pair the remaining points of $e'$. We say $e'$ is dead now. 
\end{enumerate}
Notice that in this algorithm, both the vertices and edges can have various numbers of unpaired points. Hence we must describe the expected one step changes in both the $Y_j$'s and the $Z_j$'s assuming that the minimum degree (i.e. the number of unpaired points in $v_i$ from step \tbf{I1}) is $r$. We begin with the $Y_j$'s.

\begin{align}
      &\E\sbrac{Y_j(i+1)-Y_j(i) \mid G_i \wedge \{\op_i=\Op_r \}}  \nonumber\\
& = -\bfo_{j=r} +\frac{r \cdot 2Z_2}{M_1} \rbrac{-\frac{ j Y_j}{M_1}+ \frac{ ((j+1)Y_{j+1} - jY_j)M_2L_2}{M_1^3}} + O\rbrac{\frac1{M_1}} \label{eqn:41}\\
& =   -\bfo_{j=r} - \frac{r \cdot 2Z_2 \cdot j Y_j}{M_1^2}+ \frac{r \cdot 2Z_2 ((j+1)Y_{j+1} - jY_j)M_2L_2}{M_1^4}+ O\rbrac{\frac1{M_1}}\label{eq:YExpOneStepInd}
\end{align}
The first term of line \eqref{eqn:41} accounts for the loss of a vertex in $Y_j$ if the vertex $v_i$ chosen in (\tbf{I1}) is of degree $j$ (which happens if $j=r$). The $\frac{r \cdot 2Z_2}{M_1}$ on the second term is the expected number of edges $e$ from (\tbf{I2}) which now have only one unpaired point. Each such edge triggers a vertex closure, and $\frac{jY_j}{M_1}$ is the probability that the closed vertex is in $Y_j$. Now for the last term, we see that for each such edge $e$, $\frac{M_2}{M_1} \cdot \frac{L_2}{L_1} = \frac{M_2L_2}{M_1^2}$ is the expected number of pairings in (\tbf{I4}), each of which has an expected contribution of $\frac{(j+1)Y_{j+1} - jY_j}{M_1}$. 

\begin{align}
      &\E\sbrac{Z_j(i+1)-Z_j(i) \mid G_i \wedge \{\op_i=\Op_r \}} \nn\\
      &= \frac{r((j+1)Z_{j+1}-jZ_j)}{M_1} - \frac{r \cdot 2Z_2\cdot jZ_j L_2}{M_1^3} + O\rbrac{\frac1{M_1}}\label{eq:ZExpOneStepInd}
\end{align}
The first term is the expected effect of (\tbf{I1}). The second term is the expected effect of (\tbf{I3}).

Let us emphasize here that we assume the uniformity $k$ is at least 3. We observe that this means that there will be a significant number of steps at the beginning where the minimum degree will almost always be either $d$ or $d-1$. Indeed, the degree of a vertex only decreases when one of its neighbors gets closed, which can only happen when we see an edge with $k-1$ of its vertices in the independent set. But it should take a linear number of steps before we see any significant number of such edges. Thus toward the beginning vertex closures are quite rare and the number of vertices whose degree decreases to $d-1$ is small enough that we are able to process them as they arise. We will say that \tbf{phase 0} consists of this initial string of steps when the minimum degree is always $d$ or $d-1$. Of course when phase 0 ends we will begin phase 1 and so on, where in \tbf{phase $p$} the minimum degree is almost always either $d-p$ or $d-p-1$. 

We define the scaled versions of \eqref{eq:YExpOneStepInd} and \eqref{eq:ZExpOneStepInd}:
\begin{align}
    f_{j, r}(\yb, \zb)&:= -\bfo_{j=r} - \frac{r \cdot 2z_2 \cdot j y_j}{m_1^2}+ \frac{r \cdot 2z_2 ((j+1)y_{j+1} - jy_j)m_2\l_2}{m_1^4},\\
    g_{j, r}(\yb, \zb)&:=\frac{r((j+1)z_{j+1}-jz_j)}{m_1} - \frac{r \cdot 2z_2\cdot jz_j \l_2}{m_1^3},
\end{align}
where
\begin{align}
 \l_1& =\sum_{j=1}^d jy_j, \qquad \l_2 = \sum_{j=1}^d j(j-1)y_j, \label{eqn:lind}\\
 m_1& =\sum_{j=1}^{k} jz_j, \qquad m_2 = \sum_{j=1}^{k} j(j-1)z_j.\label{eqn:mind}
\end{align}

We let
 \begin{equation}
        \alpha_p(\yb, \zb) := f_{d-p-1, d-p}(\yb, \zb), \qquad \tau_p(\yb, \zb) := -f_{d-p-1, d-p-1}(\yb, \zb).
\end{equation}
The blended differential equations for phase $p$ are then
\begin{align}
    \frac{dy_j}{dx}& = F( \yb, \zb, j, p) := 
    \frac{\tau_p}{\tau_p+\alpha_p} f_{j, d-p}( \yb, \zb) + \frac{\alpha_p}{\tau_p+\alpha_p} f_{j, d-p-1}( \yb, \zb),& \;\;\; p \le d-2\\
       \frac{dz_i}{dx}& = G( \yb, \zb, i, p) := 
    \frac{\tau_p}{\tau_p+\alpha_p} g_{i, d-p}( \yb, \zb) + \frac{\alpha_p}{\tau_p+\alpha_p} g_{i, d-p-1}( \yb, \zb),& \;\;\; p \le d-2
   \label{eqn:mixedDEQindY}
\end{align}
As in Section 3, we use the above differential equations to track several  phases of the algorithm according to what the degrees of vertices are during each phase. In particular, as before, phase $p$ ends when $\tau_p=0$. 

\section{Transferring our analysis to hypergraphs of large girth}

First, we will describe how our algorithms are equivalent to certain algorithms on bipartite graphs. 

 Let $\mc{H}=(V, E)$ be a hypergraph. We define the \tbf{bipartite graph representation of $\mc{H}$} to be the graph $G=G(\mc{H})$ with bipartition $A \cup B$ where $A=\{a_v: v \in V\}$ and $B=\{b_e: e \in E\}$. $G$ has an edge $a_vb_e$ whenever $v \in e$. If $\mc{H}$ is $k$-uniform and $d$-regular then each vertex in $A$ will have degree $d$ and each vertex in $B$ will have degree $k$, i.e. $G(\mc{H})$ will be {\bf $(d, k)$-biregular}. Note that we can generate a uniform random $(d, k)$-biregular (multi-)graph using the same pairing model from Section \ref{sec:setup}. A matching $M \subseteq E$ in $\mc{H}$ is equivalent to a set of vertices $B' \subseteq B$ all having disjoint neighborhoods. Likewise an independent set $I \subseteq V$ is equivalent to a set of vertices $A' \subseteq A$ that does not contain the neighborhood of any vertex in $B$. In fact, our two algorithms are equivalent to the following. 

 \begin{center}
\begin{minipage}{.45\textwidth}
\begin{algorithm}[H]
  \SetKwInOut{Input}{Input}
   \SetKwInOut{Output}{Output}
	\caption{\textsc{degree-greedy matching process}}
	\Input{Bipartite graph $G$ on $A \cup B$}
	\Output{$B' \subseteq B$ all having disjoint neighborhoods}
	$B' = \emptyset$\;
	\While{$B\neq \emptyset$}
	{Select $a \in A$ \uar~from all choices with minimum positive degree\;
 Select  $b \in N(a)$ \uar\;
	 $B'	\gets B' \cup \cbrac{b}$\;
	 $B\gets B\sm \{b\} \sm N^2(b)$\;
	}
	\Return $B'$\;
\end{algorithm} 
\end{minipage}
    \begin{minipage}{.45\textwidth}
\begin{algorithm}[H]
  \SetKwInOut{Input}{Input}
   \SetKwInOut{Output}{Output}
	\caption{\textsc{degree-greedy independent process}}
	\Input{Bipartite graph $G$ on $A \cup B$}
	\Output{ $A' \subseteq A$ with $N(b) \not \subseteq A'$ for all  $b \in B$}
	$A' = \emptyset$\;
	\While{$A \neq \emptyset$}
	{Select $a \in A$ \uar~from all vertices of minimum degree\;
	 $A'	\gets A' \cup \cbrac{a}$\;
	 $A \gets A\setminus \cbrac{a'\in A : \exists b\in B, N(b) \subseteq A \cup \{a'\}}$\;
	}
	\Return $A'$\;
\end{algorithm} 
\end{minipage}
\end{center}

We describe a minor adaptation of the results of Hoppen and Wormald \cite{HWComb,WPhases}, which essentially state that one can transfer the analysis of certain algorithms from random regular graphs to regular graphs of high girth. In fact, Hoppen and Wormald \cite{HWComb} noted that ``our methods can readily be extended so as to apply to $r$-regular bipartite graphs, or indeed multipartite graphs, of large girth.'' Here we just observe that their comment also applies to $(d, k)$-biregular graphs of large girth. We will sketch a proof. The setup involves only very minor changes to a few definitions, and the proofs from \cite{HWComb,WPhases} can be followed mutatis mutandis.  In the next paragraph we summarize the changes we make to the definitions. 

A \tbf{bipartite local deletion algorithm} is an algorithm that runs on a bipartite graph $G$. More specifically, the algorithm proceeds by generating a sequence of \tbf{survival graphs} $G_0=G, G_1, \ldots$, which are graphs with colors assigned to the vertices.  At each step $i$ we select a vertex $v_i$ of $G_i$, explore some vertices at distance at most $D$ from $v_i$, and form $G_{i+1}$ by possibly recoloring or deleting some of the explored vertices. We call $D$ the \tbf{depth} of the algorithm and assume it is constant. The \tbf{type} $\tau_G(v)$ of a vertex $v$ in a colored bigraph graph $G$ is a triple consisting of its color in $G$, degree in $G$ and which side of the bipartition $v$ is in. There are two sets of colors: $\mc{C}$, the \tbf{transient colors}, and $\mc{E}$, the \tbf{output colors}. The initial graph $G_0$ has all vertices the same transient color which we call \tbf{neutral}.
A \tbf{selection rule} is a function $\Pi$ that, for a
nonempty colored bigraph $G=(A \cup B,E)$, gives a probability distribution $\pi_G$ on
the power set of $A$ with the properties that
\begin{enumerate}[label=(\roman*)]
    \item For any $v, w \in A$ with $\t_G(v)=\t_G(w)$ we have
    \[
    \sum_{\substack{S \subseteq A \\ v \in S}} \pi_G(S) = \sum_{\substack{T \subseteq A \\ w \in T}} \pi_G(T),
    \]
    \item $\pi_G(S)=0$ for any $S$ containing a vertex of color $\p'$.
\end{enumerate}
Note that $\pi_G$ is a distribution on the subsets of $A$ (rather than all vertices). This is because for both of our algorithms we are selecting from only one part of the bipartition of the graph (for hypergraph matchings we select hyperedges, and for independent sets we select vertices of the hypergraph). 

The definition of a {\bf local deletion algorithm} in \cite{HWComb} uses the terminology {\em type} and {\em selection rule} but their graphs are not bipartite and so their types do not contain information about any bipartition, and for their selection rule $S$ can be any set of vertices. However, after redefining types and selection rules, we can use their definition of the term {\em local deletion algorithm}, as well as the rest of the definitions in their Section 3.1, with no further changes. 

By making the straightforward changes (just to adjust for the small changes in our definitions) to the proofs that appear in \cite{HWComb,WPhases}, we can essentially carry the analysis of our matching and independent set algorithms on random regular hypergraphs over to the setting of high girth regular hypergraphs. In the next three paragraphs we provide a high-level sketch of the proof. 

Our algorithms are {\bf prioritized} by degree, i.e. each step of our algorithm starts by selecting a vertex where the selection takes into account the current degrees of all the vertices. The first step in the proof sketch will be to {\bf deprioritize} our algorithms, meaning that we instead select vertices according to a probability distribution that depends only on time. By ensuring that this probability distribution chooses a vertex of each possible degree at approximately the same rate as the prioritized algorithm, we can derive a deprioritized version whose behavior is asymptotically the same. See Theorem 2 of \cite{WPhases}. 

Now that our algorithm is deprioritized, we {\bf chunkify} it. Our current deprioritized algorithm selects one vertex each step, but the chunkified version will select many vertices each step. We say a deprioritized local deletion algorithm is {\bf chunky} if the number of steps does not depend on the input graph. In particular this is possible because at each step $i$ the probability we select a given vertex of degree $j$ is some number $p_{i, j}$. Chunkifying an algorithm is essentially passing to a ``nibble version.'' The nibble method, famously used by R\"odl \cite{rodl} to find almost-perfect matchings in certain hypergraphs, is a technique very similar to the differential equation method where one runs a randomized process that takes large (but not too large) steps. The fact that each step is large enough for concentration results to apply sometimes makes the nibble version of a process easier to analyze. Our deprioritized algorithms can be chunkified, meaning that there exist chunky deprioritized algorithms whose performance is approximately the same as the original algorithms, analogously to Theorem 6.3 from \cite{HWComb}.

The final step of the proof is to observe that a chunky local deletion algorithm perform almost the same on all $(d, k)$-biregular bigraphs without many short cycles. This finishes the argument since random biregular bigraphs do not have many short cycles w.h.p. This part of the proof is analogous to Theorem 4.3, Lemma 4.4 and Theorem 4.5 in \cite{HWComb}.

\section{Final discussion and open questions}

First we give an alternative description of the problem we have addressed for high-girth hypergraphs. Let
\begin{align*}
\nu^*_{k, d, g} &:= \min\cbrac{\frac{\nu(\mc{H})}{|V(\mc{H})|}: \mc{H} \mbox{ is $k$-uniform, $d$-regular with girth at least $g$}},\\
\alpha^*_{k, d, g} &:= \min\cbrac{\frac{\alpha(\mc{H})}{|V(\mc{H})|}: \mc{H} \mbox{ is $k$-uniform, $d$-regular with girth at least $g$}}.
\end{align*}
It is worth noting that for any $k, d$ there exist $k$-uniform $d$-regular hypergraphs of arbitrarily high girth, and so we are never taking the minimum of the empty set above. Of course, as $g \rightarrow \infty$ we necessarily have $|V(\mc{H})|\rightarrow \infty$ as well. As $g$ increases, both $\nu_{k, d, g}$ and $\alpha_{k, d, g}$ are nondecreasing and bounded above, and so the following limits exist:
\[
\nu^*_{k, d}:= \lim_{g \rightarrow \infty} \nu^*_{k, d, g}, \qquad \alpha^*_{k, d}:= \lim_{g \rightarrow \infty} \alpha^*_{k, d, g}.
\]
Then we have shown that $\nu^*_{k, d} \ge \nu_{k, d}$ (see Table \ref{tbl1}) and $\alpha^*_{k, d} \ge \alpha_{k, d}$ (see Table \ref{tbl4}). Likewise, Tables \ref{tbl3} and \ref{tbl6} give upper bounds on $\nu^*_{k, d}, \alpha^*_{k, d}$.

Next we discuss some open questions.

\subsection{Almost-perfect matchings?}

We will restate our open question from \cite{BB}, which is about almost-perfect matchings. First we motivate the question. In the graph case $k=2$, it is known that every sequence of graphs with girth going to infinity has an almost-perfect matching. Indeed, Bohman and Frieze \cite{BF11} showed that the Karp-Sipser algorithm produces an almost-perfect matching in random regular graphs. Since Karp-Sipser is a local deletion algorithm, the results of Hoppen and Wormald \cite{HopWorm} transfer the result to high-girth regular graphs. For an alternative, non-algorithmic, approach to the same problem, see Flaxman and Hoory \cite{FH}. In the hypergraph case, Cooper, Frieze, Molloy and Reed \cite{CFMR96} proved that $\mc{H}_k(n,  d)$  has a perfect matching w.h.p.~for all $d$ large enough with respect to $k$ (and $k|n$). Therefore, we ask whether this result can be approximately carried over to all high-girth regular hypergraphs.

\begin{opq}
Fix $k \ge 3$. Does there exist some $d=d(k)$ large enough such that every sequence $\mc{H}_n$ of $k$-uniform $d$-regular hypergraphs on $n$ vertices with girth tending to infinity has an almost-perfect matching? 
\end{opq}

Since the $k=2$ case of the above question can be settled using a local deletion algorithm, it seems worthwhile to try the same approach for $k \ge 3$. Unfortunately, the degree-greedy matching process we study in this note appears to fall short of this goal (at least for the smallish values of $d$ we were able to check using Maple). However, one can be more clever and design a better algorithm. In the next paragraph we discuss one such idea (which ultimately also falls short of answering the above open question). 

We consider building upon the idea for the degree-greedy algorithm. We will still start by choosing a minimum degree vertex $v$ to get matched in this step, but we will choose the matching edge a bit more carefully. We will choose the matching edge incident with $v$ with the smallest possible sum of degrees of its vertices (since this minimizes the number of deleted edges). 
Unfortunately, this process gives rise to a system of differential equations that is significantly more complicated than the degree-greedy matching process we analyzed in this note, and the reward is quite small. We did some numerical calculations for small uniformity and degree, and the improvement does not seem worth the effort of writing more about.

\subsection{Limitations of local algorithms?}

In random regular graphs, it was conjectured by Hatami, Lov\'asz and Szegedy \cite{HLS} that local algorithms could find independent sets that are close to the maximum. Gamarnik and Sudan \cite{GS} disproved this conjecture in a result which was later improved by Rahman and Vir\'ag \cite{RV}. Indeed, \cite{RV} shows that for large $d$, there is no local algorithm which w.h.p. finds an independent set bigger than about half the size of the maximum. This fraction (half) is optimal, since even the simplest greedy algorithm achieves it. We ask the natural question about hypergraphs, preferring to keep the question a bit vague. For random regular hypergraphs with fixed uniformity and large regularity, the independence number was determined up to a small error term by the second author and Frieze \cite{BF}. 

\begin{opq}
    What are the limitations of local algorithms for independent sets (or, for that matter, matchings) in random regular hypergraphs?
\end{opq}

\subsection{Independent sets in $2$-regular hypergraphs}

Note that the $d=2$ columns of Tables \ref{tbl4} and \ref{tbl6} look like $1-\frac 1k$ (and so the lower and upper bounds match up to the third decimal). This is not too surprising, as we will now discuss. As we observed in \cite{BB}, any $k$-uniform hypergraph with $n$ vertices has independence number at most $(1-\frac 1k)n$. Let $\mc{H}$ be a high-girth $k$-uniform $d$-regular hypergraph with $d=2$. We take the {\bf dual hypergraph} $\mc{H}^*$, which has a vertex $v_e$ for each edge $e \in E(\mc{H})$ and an edge $\{v_{e_1}, \ldots v_{e_d}\}$ whenever $e_1, \ldots, e_d$ are distinct edges of $\mc{H}$ intersecting at a single vertex. Another way to think of $\mc{H}^*$ is to take the bipartite graph $G(\mc{H})$, swap the two parts of the bipartition, and take $\mc{H}^*$ to be the  corresponding hypergraph. In particular $\mc{H}^*$ is $d$-uniform and $k$-regular, so since $d=2$ then $\mc{H}^*$ is a graph. It is easy to see that $\mc{H}$ and $\mc{H}^*$ have the same girth. 

As a $k$-regular graph of high girth, $\mc{H}^*$ has an almost-perfect matching $M$. Since $\mc{H}^*$ has $2n/k$ vertices we have $|M|\approx n/k$. Furthermore, $M$ can be found using a local deletion algorithm such as Karp-Sipser. Let $V \subseteq V(\mc{H})$ be the vertex set corresponding to the complement of $M$, i.e. the edges in $\mc{H}^*$ not in $M$. Then $|V|  =  n - |M| \approx (1- \frac{1}{k}) n$. We claim that $V$ is not too far from being an independent set. Indeed, any edge $e$ of $\mc{H}$ with all $k$ of its vertices in $V$ must correspond to a vertex in $\mc{H}^*$ that is unmatched by $M$, and there are very few such vertices. For each such edge $e$ we simply remove one vertex of $e$ from $V$, obtaining an independent set $V'$ of almost the same size $\approx (1- \frac{1}{k}) n$.

To summarize: the fact that a regular graph of high girth has an almost-perfect matching (which can be found with a local deletion algorithm) implies that a 2-regular hypergraph of high girth has an independent set of size almost $(1- \frac{1}{k})n$ which can be easily found from the matching. Thus, it would not be too surprising that one can apply a local deletion algorithm directly to the 2-regular hypergraph to get such an independent set. 

It would be interesting to know whether our particular algorithm really does asymptotically achieve this upper bound. We proved that it does up to the third decimal place for $3 \le k \le 5$ (and some higher values of $k$ we did not list on the tables), but it would be nice to know precisely.  We state it as an open question. 

\begin{opq}
    Does the degree-greedy independent process produce an independent set of size $(1-\frac 1k) n + o(n)$ in $\mc{H}_k(n, 2)$?
\end{opq}

\section{Acknowledgement}

We would like to thank Carlos Hoppen and Nick Wormald for a helpful conversation.

\end{document}